\documentclass[12pt,a4paper]{article}
\usepackage{a4,amsmath,amssymb,amsfonts,amsthm}
\usepackage{amsrefs}
\usepackage[all]{xy}

\begin{document}
\bibliographystyle{plain}
 

\def\mR{\M{R}}           
\def\mZ{\M{Z}}           
\def\mN{\M{N}}           
\def\mQ{\M{Q}}       
\def\mC{\M{C}}  
\def\mG{\M{G}}



\def\Spec{{\rm Spec}}
\def\rg{{\rm rg}}
\def\Hom{{\rm Hom}}
\def\Aut{{\rm Aut}}
 \def\Tr{{\rm Tr}}
 \def\Exp{{\rm Exp}}
 \def\Gal{{\rm Gal}}
 \def\End{{\rm End}}
 \def\det{{{\rm det}}}
 \def\Td{{\rm Td}}
 \def\ch{{\rm ch}}
 \def\che{{\rm ch}_{\rm eq}}
  \def\Spec{{\rm Spec}}
\def\Id{{\rm Id}}
\def\Zar{{\rm Zar}}
\def\Supp{{\rm Supp}}
\def\eq{{\rm eq}}
\def\Ann{{\rm Ann}}
\def\LT{{\rm LT}}
\def\Pic{{\rm Pic}}
\def\rg{{\rm rg}}
\def\et{{\rm et}}
\def\sep{{\rm sep}}
\def\ppcm{{\rm ppcm}}
\def\ord{{\rm ord}}
\def\Gr{{\rm Gr}}
\def\ker{{\rm ker}}
\def\rk{{\rm rk}}
\def\Jac{{\rm Jac}}
\def\unr{{\rm unr}}
\def\rig{{\rm rig}}
\def\Fal{{\rm Fal}}


\def\beginProof{\par{\bf Proof. }}
 \def\endProof{${\qed}$\par\smallskip}
 \def\pr{^{\prime}}
 \def\prpr{^{\prime\prime}}
 \def\mtr#1{\overline{#1}}
 \def\ra{\rightarrow}
 \def\mfp{{\mathfrak p}}
 
 \def\ol#1{\overline{#1}}
 \def\mQ{{\Bbb Q}}
 \def\mR{{\Bbb R}}
 \def\mZ{{\Bbb Z}}
 \def\mC{{\Bbb C}}
 \def\mN{{\Bbb N}}
 \def\mF{{\Bbb F}}
 \def\mA{{\Bbb A}}
  \def\mG{{\Bbb G}}
 \def\mP{{\Bbb P}}  
 \def\CI{{\cal I}}
 \def\CJ{{\cal J}}
 \def\CH{{\cal H}}
 \def\CO{{\cal O}}
 \def\CA{{\cal A}}
 \def\CB{{\cal B}}
 \def\CC{{\cal C}}
 \def\CK{{\cal K}}
 \def\CL{{\cal L}}
 \def\CM{{\cal M}}
\def\CP{{\cal P}}
\def\CR{{\cal R}}
\def\CG{{\cal G}}
\def\CF{{\cal F}}
 \def\wt#1{\widetilde{#1}}
 \def\mod{{\rm mod\ }}
 \def\refeq#1{(\ref{#1})}
 \def\blb{{\big(}}
 \def\brb{{\big)}}
\def\mc{{{\mathfrak c}}}
\def\mcpr{{{\mathfrak c}'}}
\def\mcprpr{{{\mathfrak c}''}}
\def\ss{{\rm ss}}
\def\parf{{\rm parf}}
\def\P1{{{\bf P}^1}}
\def\cod{{\rm cod}}
\def\pr{\prime}
\def\prpr{\prime\prime}
\def\ss{\scriptstyle}
\def\OX{{ {\cal O}_X}}
\def\mpartial{{\mtr{\partial}}}
\def\inv{{\rm inv}}
\def\indlim{\underrightarrow{\lim}}
\def\prolim{\underleftarrow{\lim}}
\def\pprolim{'\prolim'}
\def\Pro{{\rm Pro}}
\def\Ind{{\rm Ind}}
\def\Ens{{\rm Ens}}
\def\without{\backslash}
\def\pbdb{{\Pro_b\ D^-_c}}
\def\qc{{\rm qc}}
\def\Com{{\rm Com}}
\def\an{{\rm an}}
\def\gfield{{\rm\bf k}}
\def\s{{\rm s}}
\def\dR{{\rm dR}}
\def\ari#1{\widehat{#1}}
\def\ul#1{\underline{#1}}
\def\sul#1{\underline{\scriptsize #1}}
\def\mou{{\mathfrak u}}
\def\ich{\mathfrak{ch}}
\def\cl{{\rm cl}}
\def\K{{\rm K}}
\def\R{{\rm R}}
\def\F{{\rm F}}
\def\L{{\rm L}}
\def\pgcd{{\rm pgcd}}
\def\rc{{\rm c}}
\def\N{{\rm N}}
\def\E{{\rm E}}
\def\H{{\rm H}}
\def\CHOW{{\rm CH}}
\def\A{{\rm A}}
\def\d{{\rm d}}
\def\Res{{\rm  Res}}
\def\GL{{\rm GL}}
\def\Alb{{\rm Alb}}
\def\alb{{\rm alb}}
\def\Hdg{{\rm Hdg}}
\def\Num{{\rm Num}}
\def\Irr{{\rm Irr}}
\def\Frac{{\rm Frac}}
\def\Sym{{\rm Sym}}
\def\indlim{\underrightarrow{\lim}}
\def\prolim{\underleftarrow{\lim}}
\def\Stab{{\rm Stab}}
\def\un{{\rm un}}
\def\red{{\rm red}}
\def\Per{{\rm Per}}
\def\char{{\rm char}}


\def\RHom{{\rm RHom}}
\def\rRHom{{\mathcal RHom}}
\def\rHom{{\mathcal Hom}}
\def\dotimes{{\overline{\otimes}}}
\def\Ext{{\rm Ext}}
\def\rExt{{\mathcal Ext}}
\def\Tor{{\rm Tor}}
\def\rTor{{\mathcal Tor}}
\def\SP{{\mathfrak S}}

\def\H{{\rm H}}
\def\D{{\rm D}}
\def\Del{{\mathfrak D}}

 \newtheorem{theor}{Theorem}[section]
 \newtheorem{prop}[theor]{Proposition}
 \newtheorem{propdef}[theor]{Proposition-Definition}
 \newtheorem{cor}[theor]{Corollary}
 \newtheorem{lemma}[theor]{Lemma}
 \newtheorem{sublem}[theor]{sub-lemma}
 \newtheorem{defin}[theor]{Definition}
 \newtheorem{conj}[theor]{Conjecture}

 \parindent=0pt
 \parskip=5pt

 \author{Damian R\"OSSLER\footnote{D\'epartement de Math\'ematiques, 
B\^atiment 425, 
Facult\'e des Sciences d'Orsay, 
Universit\'e Paris-Sud, 
91405 Orsay Cedex, FRANCE, E-mail: damian.rossler@math.u-psud.fr, Homepage: 
http://www.math.u-psud.fr/$\sim$rossler}}
 \title{A note on the ramification of torsion points lying on 
 curves of genus at least two}
 \date{}
\maketitle
\begin{abstract}
Let $C$ be a curve of genus $g\geqslant 2$ defined over the fraction field $K$ of 
a complete discrete valuation ring $R$ with algebraically closed residue field. Suppose 
that $\char(K)=0$ and that the characteristic of the residue field is not $2$. 
Suppose that the Jacobian $\Jac(C)$ has semi-stable reduction over $R$. 
Embed $C$ in $\Jac(C)$ using a $K$-rational point. We show 
that the coordinates of the torsion points lying on $C$ lie in the unique moderately 
ramified quadratic extension of the field generated over $K$ by the 
coordinates of the $p$-torsion points on $\Jac(C)$.
\end{abstract}

 \section{Introduction}

Let $R$ be a complete discrete valuation ring. Suppose that the residue field $k$ of 
$R$ is algebraically closed and of 
characteristic $p\geqslant 0$. Suppose that $p\not=2$. Let $K$ be the  fraction field of $R$ and suppose 
that ${\rm char}(K)=0$. 
Let $C$ be a curve of genus $g\geqslant 2$ defined over $K$. 
Let $j:C\to\Jac(C)$ be the closed immersion of $C$ into its Jacobian defined 
by a $K$-rational point. 
Let $A:=\Jac(C)$. Let $\CA$ be the N\'eron model of $A$ over $R$. 
Suppose that the connected component of the special fiber $\CA_k$ of $\CA$ is 
a semi-abelian variety (in other words, $\CA$ has semi-stable reduction). 

Let $L:=K(A[p](\bar{K}))$ be the extension of 
 $K$ generated by the coordinates of the $p$-torsion points 
 of $A(\bar{K})$. In particular, $L=K$ if $p=0$. Let $L'$ be the unique moderately ramified quadratic extension of $L$. 
 
 Finally, let $K_1\subseteq\bar{K}$ be the field generated over $K$ by the 
coordinates of the elements of $\Tor(A(\bar{K}))\cap C(\bar{K})$. Here 
$\Tor(A(\bar{K}))$ is the subgroup of $A(\bar{K})$ consisting of elements of finite order.

The aim of this note is to prove the following statement :
\begin{theor}
{\rm\bf (a)}  
$
K_1\subseteq L'
$

{\rm\bf(b)} if $C$ is not hyperelliptic then
$K_1\subseteq L$. 
\label{ThI}
\end{theor}

Theorem \ref{ThI} should be understood as a complement to 
some results of Tamagawa (see \cite{Tamagawa-Ramification}), Baker-Ribet (see \cite{Baker-Ribet-Galois}) and Coleman (see \cite{Coleman-Ramified}). 

For instance, with the present 
notation, suppose that $p>0$, that 
$R$ is the maximal unramified extension of $\mQ_p$ 
and that the abelian part of the connected component of $\CA_k$ is an ordinary abelian variety. 
Tamagawa then proves that $K_1$ is contained in the extension 
of $K$ generated by the $p$-th roots of unity (see \cite{Tamagawa-Ramification} or \cite[Th. 4.1]{Baker-Ribet-Galois}). Another example is the following result of Coleman : if $p>{\rm max}(2g,5)$, $R$ is the maximal unramified extension of $\mQ_p$ and $\CA_k$ is an abelian variety, then $K_1\subseteq K$ (see \cite[Conj. B]{Coleman-Ramified}). 

All these results restrict the size of $K_1$ under various hypotheses on 
the special fiber $\CA_k$ and on the order of absolute ramification of $K$.
The interest of Theorem \ref{ThI} is that it 
provides a limit for the size of $K_1$ under the mild hypothesis of 
semi-stability of $\CA$ only and with no assumption on 
the absolute ramification of $K$.  The hypothesis of 
semi-stability is not very restrictive, because it will automatically be satisfied, if the 
$l$-torsion points of  the Jacobian variety are $K$-rational, 
for $l$ a prime number such that $l>2$ and $l\not= p$ (Raynaud's criterion, 
see \cite[IX]{SGA7.1}). 

{\bf Remark.} The avoidance of the prime $p=2$ is critical. 
It appears in both 
Lemma \ref{lemBox} and Lemma \ref{lemWeil} and this is 
exploited at the end of the proof of Theorem \ref{ThI}. It would be interesting to 
extend this proof  to the case $p=2$. 

{\bf Notations.} If $l$ is a prime 
number and $G$ is an abelian group, we write $\Tor^l(G)$ for 
the set of elements of $\Tor(G)$ whose order is prime to $l$ and 
$\Tor_l(G)$ for the set of elements of $\Tor(G)$ whose order is a power 
of $l$. The expression $\Tor^{0}(G)$ will stand for $\Tor(G)$. 
We shall denote by $+$ the group law on $A(\bar{L})$. We shall write 
divisors on $C_{\bar{L}}$ in the form
$$
n_1 P_1\oplus n_2 P_2\oplus\dots \oplus n_r P_R
$$
where $n_i\in\mZ$. The symbol $\sim$ will be used to denote linear equivalence of 
divisors. 

\section{Proof of Theorem \ref{ThI}}

Let $L^t$ be the maximal moderately ramified 
extension of $L$. 
Let $I:=\Gal(\bar{L}|L)$, $I^w:=\Gal(\bar{L}|L^t)$ and $I^t:=\Gal(L^t | L)$. 
Recall that $I^w=0$ if $\char(k)=0$ and that there is a non-canonical isomorphism 
$I^t\simeq\oplus_{l\not=p,\ l\ {\rm prime}}\ \mZ_l$ (see \cite[chap. IV]{Serre-Local}). 
Furthermore, the group $I^w$ is a pro-$p$-group if $p>0$. 

We shall need the following five results.

\begin{theor}[monodromy theorem]
For any $x\in\Tor^p(A(\bar{L}))$ and any $\sigma\in I$, 
the equation $\sigma^2(x)-2\sigma(x)+x=0$ is satisfied.
\label{monoth}
\end{theor}
\beginProof
See \cite[IX, 5.12.2]{SGA7.1}
\endProof

\begin{lemma}
The action of $I^w$ on
$\Tor^p(A(\bar{L}))$ is trivial.
\label{lemwild}
\end{lemma}
\beginProof
This is a direct consequence of Theorem \ref{monoth}. See for instance 
\cite[Appendix, Lemma A.1]{Baker-Ribet-Galois}. 
\endProof

\begin{lemma}
The action of $I^t$ on $\Tor_p(A(L^t))$ is trivial.
\label{lemtriv}
\end{lemma}
\beginProof
We may restrict ourselves to the case where $p>0$. 
Let $T\subseteq \Tor_p(A(L^t))$ be a finite $I^t$-invariant subgroup. 
We have to show that the action of $I^t$ on $T$ is trivial.
The action of $I^t$ on $T$ preserves the order of elements, hence 
$T$ is an inner direct sum of $G^t$-invariants subgroups 
of the form $(\mZ/p^r)^s$. Hence we might suppose without loss 
of generality that $T\simeq (\mZ/p^r)^s$ for some $r,s\leqslant 1$. 
Let $T_p$ be the subgroup of $p$-torsion elements of $T$. 
The fact that the $p$-torsion points in $A(\bar{L})$ are $L$-rational 
implies that the action of $G^t$ on 
$T_p$ is trivial. Hence the image of $G^t$ in 
$\Aut(T)$ lies in the kernel of 
the natural group map
$$
\Aut(T)\to\Aut(T_p)
$$
Under the above isomorphism $T\simeq (\mZ/p^r)^s$, this corresponds 
to $s\times s$-matrices of the form 
$\Id+pM$, where $M$ is an $s\times s$-matrices with coefficients in $\mZ/p^r\mZ$. 
This last fact is a consequence of the fact that multiplication by $p^{r-1}$ induces 
an $I^t$-equivariant isomorphism $T/pT\to T_p$. 
The calculation 
$
(\Id+pM)^{p^{r-1}}=\Id
$
now shows that the image of $I^t$ in $\Aut(T)$ is 
a $p$-group. On the other hand, $I^t$ is a direct sum of pro-$l$-groups, with 
$l\not=p$. The order of the image 
of $I^t$ is thus prime to $p$. This image is thus trivial. 
\endProof

\begin{lemma}[Boxall]
Let $B$ be an abelian variety over a field $F$ of characteristic $0$. 
Let $l>2$ be a prime number and 
let $L:=F(A[l])$ be the extension of $K$ 
generated by the $l$-torsion points of $A$. 
Let $P\in \Tor_l(B(\bar{L}))$ and suppose that $P\not\in B(L)$.
Then there exists $\sigma\in\Gal(\bar{L}|L)$ 
such that $\sigma(P)-P\in B[p](\bar{L})\setminus\{0\}$.
\label{lemBox}
\end{lemma}
\beginProof
See \cite[Lemme 1]{Boxall-Sous} or \cite[Prop. 3]{Rossler-A-note}. 
\endProof

\begin{lemma}
Let $P\oplus Q$ and $P'\oplus Q'$ be two divisors of degree $2$ on
$C_{\bar{L}}$. Suppose that $P\oplus Q$ and $P'\oplus Q'$ are linearly equivalent.

If $C$ is not hyperelliptic, then the two divisors coincide.

If $C$ is hyperelliptic, then either the two divisors coincide or we have $Q=\iota(P)$ 
and $Q'=\iota(P')$.
\label{lemhyp}
\end{lemma}
Here $\iota:C\to C$ is the uniquely defined hyperelliptic 
involution.
\beginProof
See \cite[IV, Prop. 5.3]{Hartshorne-Algebraic}. 
\endProof

\begin{lemma}
Let $x\in A(\bar{L})$ and suppose that $x\not=0$. The inequality 
\begin{equation}
\#\big(C(\bar{L})\cap(C(\bar{L})+x)\big)\leqslant 2
\label{eqw1}
\end{equation} is then verified. 
If $C$ is not hyperelliptic, we even have 
\begin{equation}\#\big(C(\bar{L})\cap(C(\bar{L})+x)\big)\leqslant 1.
\label{eqw2}
\end{equation}
\label{lemWeil}
\end{lemma}
This Lemma is a consequence of \cite[Prop. 4]{Baker-Poonen-Torsion}. 
For the convenience of the reader, we provide the following proof. 
\beginProof 
We shall write $O$ for the $K$-rational point on $C$, which  
is used to embed $C$ in $A$. 


Let $a_1,\dots,a_r\in C(\bar{L})$ be pairwise distinct points 
such that $a_1+x,a_2+x,\dots,a_r+x\in C(\bar{L})$.  
Let $b_i:=a_i+x$ ($i=1,\dots,r$). 

Suppose first that $r>1$ and that 
$C$ is not hyperelliptic. 
We then have a linear equivalence
\begin{equation}
b_1\oplus a_2\sim b_2\oplus a_1
\label{lineq1}
\end{equation}
Hence either 
$b_1=b_2$ or $b_1=a_1$, either of which are ruled out. So we conclude that 
if $C$ is not hyperelliptic, then $r\leqslant 1$. This proves the inequality 
\refeq{eqw2}. 

Now suppose that $C$ is hyperelliptic and that $r>2$. Let $\iota:C\to C$ 
be the corresponding hyperelliptic involution. On top of \refeq{lineq1}, we then 
have the further linear equivalence
$$
b_2\oplus a_3\sim a_2\oplus b_3
$$
Lemma \ref{lemhyp} now implies that $a_2=\iota(b_1)$ and 
$a_2=\iota(b_3)$. Hence $b_1=b_3$, which is impossible. 
Thus we conclude that $r\leqslant 2$, if $C$ is hyperelliptic. This 
proves the first inequality \refeq{eqw1}. 
\endProof

We can now start with the proof of Theorem \ref{ThI}. 

The monodromy theorem \ref{monoth} says that 
for any $x\in\Tor^p(A(\bar{L}))$ and any $\sigma\in I$ the equation
$(\sigma-\Id)^2(x)=0$ is satisfied (remember that $\Tor^p(\cdot)=\Tor(\cdot)$ if 
$p=0$). On the other hand, Lemma \ref{lemwild} says that
$\Tor^p(A(\bar{L}))\subseteq\Tor^p(A(L^t))$ and Lemma \ref{lemtriv} implies 
that $\sigma(x)=x$ for any $x\in\Tor_p(A(L^t))$ and any $\sigma\in I^t$. Hence the equation
\begin{equation}
(\sigma-\Id)^2(x)=0
\label{sigeq}
\end{equation}
is verified for any $x\in\Tor(A(L^t))$ and any $\sigma\in I^t$. 
Let $x\in\Tor(A(L^t))\cap C(\bar{L})$ and $\sigma\in I^t$. The equation \refeq{sigeq} implies 
the linear equivalence 
\begin{equation}
\sigma^2(x)\oplus x\sim 2\sigma(x)
\label{fundeq}
\end{equation}
of divisors of degree $2$ on $C_{\bar{L}}$. 

First suppose that $C$ is not hyperelliptic; then any two linearly equivalent divisors of degree $2$ on $C_{\bar{L}}$ coincide (see Lemma \ref{lemhyp}) and the relation \refeq{fundeq} thus implies that 
$x=\sigma(x)$.\footnote{This calculation is partly the motivation for Ribet's definition 
of an "almost rational point" (see \cite[Lemma 2.7]{Baker-Ribet-Galois}) and 
is the starting point of Tamagawa's article \cite{Tamagawa-Ramification} (see Prop. 0.2 in that reference).} 

{\it Summing up, there is an inclusion}
$$
\Tor(A(L^t))\cap C(\bar{L})\subseteq A(L)
$$
{\it if $C$ is not hyperelliptic.} 

Now suppose that $C$ is hyperelliptic and let 
$\iota:C\to C$ be the uniquely defined hyperelliptic involution. 
Lemma \ref{lemhyp}  then implies that either $x=\sigma(x)$ as above or 
that the equations $\iota(\sigma(x))=\sigma(x)$ and $\iota(x)=\sigma^2(x)$ hold. Suppose 
the latter. Since $\iota$ is defined 
over $L$, the equation $\iota(\sigma(x))=\sigma(x)$ implies that 
$\iota(x)=x$. This together with the equation $\iota(x)=\sigma^2(x)$ implies that 
$\sigma^2(x)=x$. 

Now notice that since the group $I^t$ is abelian, the set  
$(I^t)^2$ of squares of elements is a normal subgroup of $I^t$. Let 
$J'$ be the Galois extension of $L$ defined by $(I^t)^2$. 
Since \mbox{$I^t\simeq\oplus_{l\not=p,\ l\ {\rm prime}}\ \mZ_l$,} we see that 
$[J':L]=2$. In the last paragraph, we showed that 
$\Tor(A(L^t))\cap C(\bar{L})$ is fixed by $(I^t)^2$. 
In other words, we have shown that $\Tor(A(L^t))\cap C(\bar{L})\subseteq A(J')$ 
where $J'$ is the unique moderately ramified quadratic extension of $L$. 
In the notation of the introduction, $J'=L'$. 

{\it Summing up, we see that}
$$
\Tor(A(L^t))\cap C(\bar{L})\subseteq A(L')
$$
{\it if $C$ is hyperelliptic.}

If $p=0$, then $L^t=\bar{K}$ so this completes the proof of Theorem \ref{ThI} in that case.

Now let $x\in\Tor(A(\bar{L}))\cap C(\bar{L})\backslash C(L^t)$. 
Let $x=x^p+x_p$ be the decomposition of $x$ into its 
components of prime-to-$p$ torsion and $p$-primary torsion, respectively. 
Lemma \ref{lemtriv} implies that $x_p\in A(\bar{L}))\backslash A({L}^t)$. 
Also, using Boxall's lemma \ref{lemBox} and the fact that $L$ contains the coordinates of 
the $p$-torsion points, we see that there exists 
$\sigma_x\in I^w$ such that $\sigma_x(x)-x=\sigma_x(x_p)-x_p\in A[p](L)\backslash\{0\}$. 
Hence 
\begin{equation}
\sigma_x(x)\in \bigcup_{\tau\in A[p]({L})\backslash\{0\}}C(\bar{L})\cap \big(C(\bar{L})+\tau\big).
\label{basic1}
\end{equation}
Lemma \ref{lemWeil} now implies that 
\begin{equation}
\#(C(\bar{L})\cap \big(C(\bar{L})+\tau\big))\leqslant 2
\label{basic2}
\end{equation}
for all $\tau \in A[p](L)\backslash\{0\}$. 
Notice that if $\sigma_x(x)\in C(\bar{L})\cap \big(C(\bar{L})+\tau_0\big)$ 
for a particular $\tau_0 \in A[p](L)\backslash\{0\}$, then we have 
$$
\{\sigma_x(x),\sigma^2_x(x),\dots,\sigma^{p}_x(x)\}\subseteq C(\bar{L})\cap \big(C(\bar{L})+\tau_0\big)
$$
(remember that by construction $\tau_0$ is fixed by $I^w$). On the other hand
$$
\{\sigma_x(x),\sigma^2_x(x),\dots,\sigma^{p}_x(x)\}=\{x+\tau_0,x+2\tau_0,\dots,x+p\tau_0=x\}.
$$
Since $p>2$ and $\tau_0$ has exact order $p$ in $A(\bar{K})$, this leads to a contradiction. 
Thus we have 
$$
 \Tor(A(\bar{L}))\cap C(\bar{L})\subseteq A(L^t).
 $$
Now remember that we have shown above (see the italicized sentences) that 
{$\Tor(A({L}^t))\cap C(\bar{L})\subseteq A(L)$} if $C$ is not 
hyperelliptic and that $\Tor(A({L}^t))\cap C(\bar{L})\subseteq A(L')$ 
if $C$ is hyperelliptic. This concludes the proof of (a) and (b).

\end{document}